\documentclass{scrartcl}

\usepackage[T1]{fontenc}
\usepackage[latin1]{inputenc}
\usepackage{amssymb,amsmath}

\newtheorem{thm}{Theorem}[section]
\newtheorem{cor}[thm]{Corollary}

\newtheorem{prop}[thm]{Proposition}
\newtheorem{defn}[thm]{Definition}
\newtheorem{rem}[thm]{Remark}

\newcommand{\abs}[1]{\left\vert#1\right\vert}

\newcommand{\Real}{\mathbb R}
\newcommand{\Cp}{\mathbb C}
\newcommand{\eps}{\varepsilon}

\newcommand{\C}{\mathcal{C}}

\newcommand{\del}{\partial}
\newcommand{\db}{\overline\partial}

\renewcommand{\rho}{\varrho}
\renewcommand{\phi}{\varphi}
\def\carre{\hbox{
\vrule height 1.453ex  width 0.093ex  depth 0ex
\vrule height 1.5ex  width 1.3ex  depth -1.407ex\kern-0.1ex
\vrule height 1.453ex  width 0.093ex  depth 0ex\kern-1.35ex
\vrule height 0.093ex  width 1.3ex  depth 0ex}\,}
\begin{document}

\title {On the displacement rigidity of Levi flat hypersurfaces - The case of boundaries
of disc bundles over compact Riemann surfaces 
}%
\author{Klas Diederich \and Takeo Ohsawa}%
%\address{}%
%\email{}%

%\thanks{}%
%\subjclass{}%
%\keywords{}%

\date{April 15, 2005 }%
%\dedicatory{}%
%\commby{}%
% ----------------------------------------------------------------

\maketitle
\begin{abstract}
{\noindent\sffamily\bfseries\large Abstract}\\
Non-existence theorems for Levi flat hypersurfaces have found
great interest in the literature. The question next to this that
has to be asked is, when existing Levi flat hypersurfaces are at
least rigid under deformations. Here, the case of boundaries of
disc bundles over certain compact Riemann surfaces is considered.
\end{abstract}
% ----------------------------------------------------------------
\section*{Introduction}

In several complex variables, Levi flat hypersurfaces arose as counterexamples in generalized
function theory on complex manifolds (cf. \cite{N3,O7,Ba2,DF3}). Nowadays, they are considered as objects
of independent interest, since they also arise as typical examples of minimal closed subsets
consisting of leaves of complex analytic foliations (cf. \cite{Ce1}). However, very few Levi
flat hypersurfaces have yet been analyzed. The most remarkable results are nonexistence
theorems (cf. \cite{S11},\cite{S12},\cite{Iord1}). Recently an attempt has begun to classify them
(cf. \cite{MatO}, \cite{O14}).
Under these circumstances we would
like to continue the study of disc bundles over compact Kähler manifolds from \cite{DO3}, where
we proved as main result that any holomorphic disc bundle $D$ over a compact Kähler manifold
$M$ is weakly $1$-complete (i.e. it admits a $\mathcal{C}^\infty$-plurisubharmonic exhaustion
function). \smallskip \\
What we want to pursue further, is a rigidity property of $\partial D$ when $D$ is
identified with a domain in the associated $\mathbb{P}^1$-bundle, say $P \rightarrow M$. For
that we shall restrict ourselves here, as a first step, to the case where $M$ is a compact Riemann surface
of genus $g > 1$, since we can exploit some (deep) results on Riemann surfaces. \smallskip \\
First, employing
Schoen-Yau's diffeomorphism theorem for harmonic maps, we refine the previous result as follows.
\begin{prop}[Consequence of Proposition \ref{Takeuchiconvexity}]
Let $C$ be a compact Riemann surface of genus $g>1$
and let $D \rightarrow C$ be a holomorphic disc bundle associated to a homomorphism $\varrho$
from the fundamental group of $C$ into the automorphism group of the unit disc $\mathbb{D}$. If
the image $\Gamma$ of $\varrho$ is a Fuchsian group such that $\mathbb{D}/\Gamma$ is homoeomorphic
to $C$, then $D$ is Takeuchi $1$-convex in the associated $\mathbb{P}^1$-bundle $P$. (For the
definitions see $\S 2$)
\end{prop}

Based on this observation, we can conclude the following rigidity result.
\begin{thm}
$\partial D$ is rigid in $P$ if either $\Gamma$ is an abelian
group or a Fuchsian group such that $\mathbb{D}/\Gamma$ is
biholomorphic to $C$ or to its conjugate $\overline C$.
\end{thm}
(For the definition of rigidity, see $\S 2$.) For the proof of the
theorem we need a Hartogs type extension theorem of Ivashkovitch
\cite{Iv1} and a basic fact on projective structures described by
R. C. Gunning in \cite{G3,G4}.

The condition on $\Gamma$ does not seem to be really essential for the rigidity of $\partial D$.
One might even suspect that the rigidity holds true for any disc bundle over any compact complex
manifold. Although there are few methods available to study the question in such generality,
the authors believe that the present work gives some insight towards that direction.

\underline{Acknowledgement}: The second-named author would like to
express his gratitude to the university of Wuppertal for the
hospitality during the preparation of this work. The authors are
very grateful to the referee for communicating to them equality
(6) which they were not aware of.

\section{Disc bundles over compact Kähler manifolds - review and refinement}

Let $M$ be a compact complex manifold of dimension $n$ and let $D \rightarrow M$ be a holomorphic fiber
bundle the fibers of which are biholomorphic to the unit disc $\mathbb{D}:=\{\zeta \in \mathbb{C}: |\zeta|<1\}$.
Recall that the group $\mathrm{Aut}\,\mathbb{D}$ of biholomorphic automorphisms of $\mathbb{D}$ consists
of the maps
$$\zeta \mapsto e^{i\theta}\frac{\zeta - a}{{\overline a} \zeta -1}\qquad (\theta \in \mathbb{R},\; a \in \mathbb{D})$$
so that transition maps of the bundle $D$ are locally constant. Accordingly, the pull-back of $D$ to the
universal covering $\widetilde{M} \rightarrow M$ is the trivial disc bundle $\tilde M \times \mathbb{D}$,
and there exists a homeomorphism
\begin{equation} \label{rho}
\varrho : \pi _1(M,x_0) \rightarrow \mathrm{Aut}\,\mathbb{D} \qquad (x_0 \in M)
\end{equation}
uniquely determined up to the inner automorphism operation of $\mathrm{Aut}\,\mathbb{D}$, such that
\begin{equation} \label{quotient}
D \simeq \widetilde{M} \times \mathbb{D}/\sim
\end{equation}
where $(x,\zeta) \sim (y,\zeta '):\Leftrightarrow y=\sigma (x)$ and ${\zeta '}= \varrho (\sigma ) \zeta$
for some $\sigma \in \pi _{1}(M,x_0 )$. Here the action of $\sigma$ is defined as the covering transformation.
We shall denote $D$ by $D_\rho$ when we want to refer to $\rho$. In particular, $D$ is uniquely determined
by $\rho$. \smallskip \\
Let $P \rightarrow M$ be the $\mathbb{P}^1$-bundle associated to the composite of $\rho$ and the inclusion
homomorphism $\mathrm{Aut}\,\mathbb{D} \hookrightarrow \mathrm{Aut}\,\mathbb{P}^1$. $D$ will be naturally
identified with a domain in $P$.
\begin{prop} \label{cases}
For any compact Kähler manifold $M$ and for any disc bundle $D$ over $M$, one of the following four cases
occurs:
\begin{enumerate}
\item $D$ admits a unique locally nonconstant pluriharmonic section.
\item $D$ admits a locally constant section.
\item $\partial D \rightarrow M$, the associated circle bundle, admits a unique
locally constant section.
\item $\partial D \rightarrow M$ admits precisely two locally constant sections.
\end{enumerate}
\end{prop}
{\em Proof:} It is clear that 2) occurs if and only if $D$ is the tubular neighborhood
of the zero section of a topologically trivial line bundle over $M$. 3) (resp. 4)) occurs
if and only if the transition maps have one (resp. two) common fixed points on the
boundary of $\mathbb{D}$, in which case there are no (resp. infinitely many)
pluriharmonic sections of $D$. By \cite{DO3} the rest is contained in the (possibly empty) case
1). \hfill\carre
\begin{cor} [cf. \cite{DO3}]
Any holomorphic disc bundle over a compact Kähler manifold is weakly $1$-complete.
\end{cor}
\begin{defn}
A $\mathcal{C}^2$ real valued function $\varphi$ on a complex manifold $X$ of
dimension $n$ is said to be $q$-convex at a point $x \in X$ if the Levi form
of $\varphi$ has at least $n-q+1$ positive eigenvalues at $x$.
\end{defn}
In the cases 1), 3) and 4) of Proposition \ref{cases} for the disc bundles $D$ over a compact Kähler
manifold of dimension $n$, it turns out that $D$ admits an exhaustion function of class
$\mathcal{C}^\infty$ which is $n$-convex outside a compact subset of $D$. This can easily be seen
from the proof of the above Corollary which we have given in \cite{DO3}. \\
Instead of this general fact, we shall prove a refined variant which we shall need later.
\begin{defn}
A relatively compact domain $\Omega$ with $\mathcal{C}^2$-smooth boundary in a complex manifold
$X$ is said to be Takeuchi $q$-convex if $\Omega$ admits a defining function $r$ of class
$\C ^2$ such that, with respect to some Hermitian metric on $X$, at least $n-q+1$ eigenvalues
of the Levi form of $-\log (-r)$ are greater than $1$ outside a compact subset of $\Omega$.
\end{defn}
\begin{rem}
A. Takeuchi \cite{Ta1} was the first to verify that the $q$-convexity in the above sense holds
for $q=1$, if $X=\mathbb{P}^n$ and $\Omega$ is a proper locally pseudoconvex domain.
\end{rem}
\begin{prop} \label{Takeuchiconvexity}
Let $C$ be a compact Riemann surface of genus $\geq 2$, and let $D
\rightarrow C$ be a disc bundle with a harmonic section $s:C
\rightarrow D$. Suppose that the set of critical points of $s$ is
finite. Then $D$ admits a defining function $r$ in the associated
$\mathbb{P}^1$-bundle $P$ such that $\del\db (-\log (-r))$
dominates the ambient metric near $\del D$. In other words, $D$ is
Takeuchi $1$-convex in $P$.
\end{prop}
{\em Proof:} In the above situation we define a function $\phi$ on $D$ by putting
$$\phi (x,\zeta):=-\log \left(1-\abs{\frac{\zeta -h(x)}{\overline{h(x)}\zeta -1}}^2 \right)$$
in terms of the coordinates $\zeta$ on the fibers. \\
Clearly, the value $\phi$ does not depend on the choice of the coordinates $\zeta$,
and $\phi$ is a real analytic exhaustion function of $D$. \\
Let $x_0$ be any point of $C$ and let $\zeta$ be a coordinate on the fiber $D_{x_0}$ satisfying
$h(x_0)=0$. Then a simple calculation gives with respect to a local coordinate $z$ around $x_0$:
\begin{eqnarray}
  \phi _{z\overline z}&=&\abs{h_z}^2 +\abs{h_{\overline z}}^2 -2 \mathrm{Re}\, (\overline\zeta ^2 h_z h_{\overline z}) \\
  \phi _{\zeta\overline\zeta}&=&\left( 1-\abs{\zeta}^2 \right)^{-2} \label{phizeta}\label{zeta}\\
  \phi _{z\overline\zeta}&=&-h_z
%\end{array}
\end{eqnarray}
at $(x_0,\zeta)$. Hence \medskip \\
\begin{equation}\label{neu}
\phi_{z\overline z}\phi_{\zeta\overline\zeta}-\abs{\phi
_{z\overline\zeta}}^2=\left( 1-\abs{\zeta}^2\right) ^{-2} \left[
\abs{\zeta}^2 (1-\abs{\zeta}^2 ) \abs{h_z}^2 +\abs{\overline\zeta
h_z-\overline{h_{\overline z}}}^2\right]
\end{equation}
holds true.\medskip \\
Let $\{x_1,\ldots ,x_m\}$ be the set of critical points of $s$,
and let $\zeta _i$ be fiber coordinates of $D$ over neighborhoods
$U_i$ of $x_i$.\\
Let $\rho_i$ be a $\mathcal{C}^{\infty}$ nonnegative functions on
$C$ whose supports are contained in $U_i$, such that $\rho _i$ is
identically $1$ near $x_i$.\\
We put, for $\varepsilon >0$,
$$\Phi=\phi +\varepsilon\sum_{i=1}^m\left\{-\rho _i\log (1-\abs{\zeta _i}^2)+
(1-\rho _i)\phi\right\}$$ Then (\ref{zeta}) and (\ref{neu}) imply
that there exists a Hermitian metric $g_P$ on the associated
bundle $P\rightarrow C$ such that $\del\db\Phi >g_P|D$ holds
outside a compact subset of $D$, if $\eps$ is chosen to be
sufficiently small.
\begin{rem} \label{refl}
Notice that for a disc bundle $D\rightarrow C$ associated to the homomorphism
$\rho :\pi _1 (C,x_0)\rightarrow \mathrm{Aut}\, \mathbb{D}$, the exterior of $D$ in the
associated $\mathbb{P}^1$-bundle $P\rightarrow C$ is equivalent to the disc bundle
associated to the homomorphism $\overline\rho$ defined by
$$\overline\rho (\sigma)(z):= \overline{\rho (\sigma )(\overline z)}$$
Therefore, if $(x,h(x))$ is a locally quasiconformal section of $D$, $(x,\overline{h(x)})$ is
a locally antiquasiconformal section of $D_{\overline\rho}$.
\end{rem}
In order to deduce the Proposition of the introduction from Proposition \ref{Takeuchiconvexity}, we
note that a $\C ^\infty$-section of $D_\rho \rightarrow C$ is naturally identified with a
$\C ^\infty$-map from $C$ to $\mathbb{D}/\mathrm{Im}\, \rho$ which is homotopic to a diffeomorphism.
Therefore, the harmonic section of $D_\rho\rightarrow C$, which exists according to the theorem
of Eells and Sampson (cf. \cite{ES}), is either quasiconformal or antiquasiconformal as
a map to $\mathbb{D}/\mathrm{Im}\,\rho$,
since it is a diffeomorphism in virtue of a theorem of Schoen-Yau \cite{SY}. Hence, the required
Takeuchi $1$-convexity follows from Proposition \ref{Takeuchiconvexity}.
\section{Stability of $q$-convexity}
Before starting to discuss the rigidity property of Levi-flat hypersurfaces, we shall prove the stability
of Takeuchi $q$-convexity for domains with Levi-flat boundaries. \smallskip \\
Recall that a $\C ^2$-smooth real hypersurface $S$ in a complex manifold $X$ is said to be Levi-flat
if $S$ locally admits a defining function the Levi form of which, restricted to the holomorphic
tangent space of $S$, is identically $0$. A real hypersurface of class $\C ^\omega$ then is defined
by pluriharmonic functions if and only if it is Levi-flat.
\begin{prop} \label{Takinvar}
Let $\Omega \subset\subset X$ be a Takeuchi $q$-convex domain with a Levi-flat $\C ^\omega$ hypersurface
as boundary. Then any $\C^2$ small perturbation of $\Omega$ as a domain with $\C ^\omega$ Levi-flat boundary
also is Takeuchi $q$-convex. In other words, there exists a tubular neighborhood $U$ and a $\C^\omega$
diffeomorphism $\phi$ between $U$ and the normal bundle of $\del\Omega$ which identifies $\del\Omega$
with the zero section of the bundle,
such that a $\C ^\omega$ domain
$\Omega ' \subset\subset X $ is Takeuchi $q$-convex if $\del\Omega '$ is Levi flat and sufficiently small
as a $\C^2$ section of the bundle.
\end{prop}
{\em Proof:} Let $r$ be a $\C ^2$ defining function of $\del\Omega$ such that $\del\db (-\log (-r))$
has at least $n-q+1$ eigenvalues $>1$ near $\del\Omega$ with respect to a Hermitian metric on $X$.
Since $\del\Omega$ is $\C ^\omega$ and Levi flat, $r$ is locally the product of a pluriharmonic defining
function, say $r_\alpha$, and a positive $\C ^2$ function, say $u_\alpha$. If $\Omega '$ is a domain with
$\C ^\omega$ Levi flat boundary such that $\del\Omega '$ is sufficiently close to $\del\Omega$ in the
$\C ^\omega$ topology, then one can choose locally pluriharmonic defining functions of $\del\Omega '$ which
are close to $r_\alpha$ even in the $\C ^\omega$ topology. Therefore,
$\Omega '$ admits a defining function $r'$ whose pluriharmonic and positive factors are close to the
corresponding factors $r_\alpha$ and $u_\alpha$ in the $\C ^2$ topology, respectively. \\
Combining this observation with
\begin{eqnarray*}
\del\db (-\log (-r))&=&\del\db (-\log (-r_\alpha u_\alpha)) \\
 &=&\del\db (-\log (-r_\alpha)) + \del\db (-\log u_\alpha)) \\
 &=&\frac{\del r_\alpha\db r_\alpha}{r^2_\alpha}+\del\db (-\log u_\alpha)
\end{eqnarray*}
one can immediately see the validity of the conclusion because of the continuity of the
eigenvalues of $\del\db (-\log u_{\alpha})$. \hfill \carre \medskip \\
\begin{defn} \label{displ}
In what follows we say that $\del\Omega '$ is a displacement of
$\del\Omega$ if $\del\Omega '$ is identifiable with a section of
the normal bundle of $\del\Omega$ by means of some diffeomorphism
$\phi$ as in Prop. \ref{Takinvar}. "Sufficiently small" always
refers to the $\C^2$-norm. $\del\Omega$ is said to be rigid if,
for any fixed choice of $\phi$, any sufficiently small
displacement of $\del\Omega$ is isomorphic to $\del\Omega$ as a
$\mathrm{CR}$ manifold.
\end{defn}
\begin{prop} \label{Takinv1}
Let $D\rightarrow C$ be a disc bundle as in Proposition \ref{Takeuchiconvexity}. Then any small
$\C^\omega$ Levi flat displacement of $\del D$ in the associated $\mathbb{P}^1$ bundle $P\rightarrow C$
bounds a Takeuchi 1-convex domain.
\end{prop}
\section{Proof of the rigidity theorem}
First we will prove case 1) of the Theorem. \medskip \\
Let $C$ be a compact Riemann surface of genus $g \geq 2$ and let $D_\rho \rightarrow C$ be a disc bundle
such that the image $\Gamma$ of $\rho$ satisfies $\mathbb{D}/\Gamma \simeq C$ or
$\mathbb{D}/\Gamma \simeq \tilde{C}$. By remark \ref{refl} it suffices to prove the theorem in the case
$\mathbb{D}/\Gamma \simeq C$, which we shall now consider.\\
Since $D_\rho$ is according to the assumption biholomorphic to the quotient
of $\mathbb{D}\times\mathbb{D}$ modulo the diagonal action
of $\Gamma$ by $(z,z)\mapsto (\gamma (z),\gamma (z))$ for $\gamma\in\Gamma$, $D_\rho$ admits a holomorphic
section $s:C\rightarrow D_\rho$ corresponding to the diagonal $\Delta =\{(z,z):z\in\mathbb{D}\}$. \\
Suppose now there exists a sequence of real analytic Levi flat hypersurfaces $S_k$, $k=1,\ldots $, converging to
$\del D_\rho$ in the $\C^2$ sense. Then by Proposition \ref{Takinv1}, $S_k$ ($k\gg 1$) separates the associated $\mathbb{P}^1$
bundle $P_\rho$ into two Takeuchi $1$-convex domains, say $D^{+}_k$ and $D^{-}_k$. We normalize notation
such that $D^{+}_k \supset s(C)$.\smallskip \\
We note that the domain $D^{-}_k$ does not contain any compact complex curve.
In fact, if there were such a curve
$C_k \subset D^{-}_k$, $C_k$ would define a multivalued holomorphic section of the affine line bundle
$P_\rho \setminus s(C)$. By averaging $C_k$ produces a holomorphic section, so that the bundle becomes
a line bundle. On the other hand, it is known that $\rho$ lifts to a $\mathrm{GL}(2,\Cp )$ representation,
say $\tilde\rho$, and the rank two vector bundle $V$ associated to $\tilde\rho$ admits a flat connection
(cf. \cite{G3}). However, $s$ and the zero section of $P_\rho \setminus s(C)$ lift to line subbundles
$L_\infty$ and $L_0$ of $V$, so that $V$ is holomorphically equivalent to $L_0 \oplus L_\infty$.
This means that, by Weil's criterion on the existence of flat connections (cf. \cite{G4}), $L_0$ and $L_\infty$
both admit flat connections, which is an absurdity because $\abs{\mathrm{deg}\, L_0}=\abs{\mathrm{deg}\, L_\infty}=
2g-2\neq 0$. Hence $D^{-}_k$ is Stein. \smallskip \\
On the other hand, since $S_k$ is $C^{\omega}$ Levi flat, there exists a neighborhood $U_k \supset S_k$ and a
complex analytic foliation $\mathcal{F}_k$ on $U_k$ of codimension $1$ which extends the foliation on $S_k$
defined by the holomorphic tangent bundle of $S_k$. We note that $\mathcal{F}_k$ is naturally identified
with a holomorphic map from $U_k$ to the projectivization of the tangent bundle of $P$. Hence, in virtue
of the extension theorem of Ivashkovitch (cf. \cite{Iv1}), $D^{-}$ admits a holomorphic foliation, possibly
with finitely many singularities, which extends $\mathcal{F}_k$. Similarly, since $D^{+}_k$ does not contain
any compact complex curves other than $s(C)$ ($s(C)$ would not be exceptional otherwise), $\mathcal{F}_k$
extends to a holomorphic foliation on $D^{+}_k \setminus s(C)$, possibly with finitely many singularities.
But since $s(C)$ is exceptional, the foliation further extends, in virtue of the Remmert-Stein continuation
theorem for complex analytic subsets (cf. \cite{S13}), to $D^{+}_k$, possibly with finitely many singularities.
Since these foliations converge to a foliation $\mathcal{F}$ of $P_\rho$ consisting of locally flat sections,
they have for sufficiently large $k$ no singularities. (Any small perturbation of a holomorphic
section as a meromorphic section is holomorphic.)\\
Thus we obtain a family of holomorphic foliations, say $\tilde{\mathcal{F}}_k$ on $P_\rho$, for $k\gg 1$,
extending $\mathcal{F}_k$ and converging to $\mathcal{F}$. Note that $\mathcal{F}_k$ defines a flat structure
on $P_\rho$. Let $\rho _k$ be the corresponding representation of $\pi _1 (C,x_{0})$ into $\mathrm{PSL}(2,\Cp )$.
Then $\del D^{+}_k \cap P_{x_0}$ must be a circle in the Riemann sphere $P_{x_0}$ because it is a closed simple
closed curve which contains an orbit of a point, consisting of infinitely many points, through the action of
$\gamma ^l$ ($l=1,2,\ldots )$ for some $\gamma \in \mathrm{PSL}(2,\Cp)$. Therefore, $D^{\pm}_k$ are
biholomorphically equivalent to some disc bundles over $C$.\\
This means, that the $\rho _k$ are equivalent to $\mathrm{Aut}(\mathbb{D})$-representations of
$\pi _1 (C,x_{0})$, say $\rho '_k$. Then $\rho '_k$ must be $\mathrm{Aut}(\mathbb{D})$-equivalent
to $\rho$, because $D_{\rho '_k}$ and $D_\rho$ both contain $s(C)$ as a holomorphic section. Therefore,
there exist bundle equivalences $\phi _{k}:D_\rho \rightarrow D_{\rho '_k}$ which extend, by inversion,
to bundle automorphisms $\tilde{\phi}_k$ of $P_\rho$ ($=P_{\rho _k} =P_{\rho '_k}$).\medskip \\
Now we are going to prove case 2) of the Theorem. Since $\Gamma$ now is abelian, the transition maps of
$D$ are either all elliptic, all parabolic or all hyperbolic (if they are not the identity).\\
Suppose first that either $D$ is trivial or the transition maps are all elliptic, and let $\{\zeta _\alpha \}$
be a system of fiber coordinates of $D$ subordinate to an open covering $\{U_\alpha \}$ of $C$ such that
$\zeta _\alpha =e^{i\theta _{\alpha\beta}}$ for $\theta _{\alpha\beta} \in \Real$ over $U_\alpha \cap U_\beta$.\\
Let $\eps >0$ and let $S$ be any $\C ^\omega$ Levi flat displacement of $\del D$ in the domain
$\{1-\eps <|\zeta_\alpha | < 1+\eps \}$. Let $\Omega _+$ and $\Omega _-$ be the components of $P\setminus S$,
such that $\Omega _+ \subset \{|\zeta_\alpha |<\infty\}$. We define a continuous function $\delta$ on
$\{|\zeta_\alpha |<\infty\}$ by letting
\begin{equation}
\delta (x,\zeta _\alpha)=\inf \{|\zeta _\alpha -\zeta '_\alpha |:(x,\zeta '_\alpha )\in S\}
\end{equation}
Then $-\log\delta + |\zeta _\alpha |^2$ is plurisubharmonic on $\Omega _+$ and real analytic near $S$.
Since $\Omega _+$ is clearly not 1-convex, $-\log\delta$ must depend only on the fiber coordinate
(cf. \cite{DO2}). Hence $S$ is CR-equivalent to $\del D$ by the fiberwise retraction along the radial
directions with respect to $\zeta_\alpha$.\smallskip \\
Next, suppose that the fiber maps are all parabolic, and let $\sigma :C\rightarrow \del D$ be the locally
constant section consisting of the common fixed points of the transition maps. In view of the classification
in Proposition \ref{cases} and the extension argument for $D_\rho$ with Fuchsian representation $\rho$,
it suffices to show that any small $\C ^\omega$ Levi flat displacements of $\del D$ are the boundaries
of 1-convex domains from both sides. In fact, if $D'\rightarrow C$ is another disc bundle such that the
associated $\mathbb{P}^1$ bundle $P'\rightarrow C$ is biholomorphic to $P$, $P$ and $P'$ are equivalent
$\mathbb{P}^1$ bundles because the genus of $C$ is not zero. Hence, the transition maps of $D'$ are all
parabolic, too. Hence $\del D'$ must contain $\sigma (C)$, so that there exists a biholomorphism between
$D$ and $D'$ given by the fiberwise translations.\\
Let $\Omega _+$ and $\Omega _-$ be the connected components of $P\setminus S$ which we want to prove
to be 1-convex. If $S\supset \sigma (C)$, there is nothing left to prove because $P\setminus \sigma (C)$
is already Stein. So let us suppose $S\not\supset \sigma (C)$. Then $S\cap\sigma (C)$ is either empty
or the union of finitely many irreducible real analytic curves, say $\gamma_1,\ldots,\gamma_N$, because of
the real analyticity of $S$. If $S\cap\sigma(C)=\emptyset$, $\Omega_+\supset\sigma (C)$ or
$\Omega_-\supset\sigma (C)$. In any case, since $\Omega_\pm\setminus\sigma(C)$ are Stein, there would arise a
2-dimensional Stein manifold with a disconnected boundary which is absurd.\\
Hence $S\cap\sigma(C)\neq\emptyset$. Let $\pi$ be the bundle projection $P\rightarrow C$. Then $\pi(\gamma_k)$
are all real analytic curves, so that $C\setminus\bigcup\pi(\gamma_k)$ carries a bounded strictly subharmonic
function say $\psi$.\\
On the other hand, let $\zeta_\alpha$ be the fiber coordinates of $P\setminus\sigma(C)$ whose transition
relations are $\zeta_\alpha=\zeta_\beta+\xi_{\alpha\beta}$, $\xi_{\alpha\beta}\in\Cp$. Hence we define
a function $\delta$ on $\Omega_\pm \setminus\sigma(C)$ by
$\delta (x,\zeta_\alpha):=\inf \{\, |\zeta_\alpha -\zeta '_\alpha|:(x,\zeta '_\alpha )\in S\}$
and obtain a plurisubharmonic function $\phi$ on $\Omega _\pm$ by extending $-\log\delta$
to $\sigma(C)\setminus S$ as $-\infty$.\\
Let $d(p)$ be the distance from $p$ to $S$ with respect to some real analytic metric $g$ on $P$,
let $\chi:P\rightarrow [0,1]$ be a $\C^\infty$ function such that $\chi\equiv 1$ on a neighborhood
of $S\cap\sigma(C)$ and that there exists a strictly plurisubharmonic function, say $\eta$, on a neighborhood
of $\mathrm{supp}\,\chi$. Then we put
\begin{equation}
\Phi:=\max\left\{\lambda(-\log\delta+\psi),-\log d+\chi\eta \right\}
\end{equation}
for a $\C^\infty$ convex increasing function $\lambda$ with $\inf\lambda=0$.\\
Clearly, if $S$ is a sufficiently small displacement of $\del D$, $\Phi$ is an exhaustion
function of $\Omega_\pm$ and satisfies
\begin{equation}
\del\db\Phi > cg
\end{equation}
near $S$ for some positive constant $c$, in the distribution sense. Therefore $\Omega_\pm$ are
1-convex.\smallskip\\
Finally, suppose that the transition maps are hyperbolic. Then $P\rightarrow C$ admits two holomorphic
curves $C_1$ and $C_2$, lying in $\del D$, which consist of the common fixed points of the transition maps.
Then the $\Cp ^*$ bundle $P\setminus (C_1 \cup C_2)$ admits a system of fiber coordinates $\{\zeta_\alpha\}$
with transition relations $\zeta_\alpha =e^{i\theta_{\alpha\beta}}\zeta_\beta$ ($\theta_{\alpha\beta}\in\Real$).
Hence, for any sufficiently small $\C ^\omega$ Levi flat displacement $S$ of $\del D$, the components of
$P\setminus S$ are 1-convex, similar as in the parabolic case. The rest also is similar to that case.\hfill Q.E.D.
\makeatletter \renewcommand{\@biblabel}[1]{\hfill#1.}\makeatother
\newcommand{\bysame}{\leavevmode\hbox to3em{\hrulefill}\,}

% ----------------------------------------------------------------
%\bibliographystyle{C:/bibtex/olko}
%\bibliography{E:/Arbeit/woc/pretex}
\renewcommand{\baselinestretch}{0.7} \Large \normalsize \ \par
  \vspace{-\baselineskip}\begin{tabbing}
 \=\mit \footnotesize Klas\ Diederich\ \ \ \ \ \ \ \ \ \ \ \ \ \ \ \ \ \ \ \ \ \ \ \ \ \ \ \ \ \ \ \ \ \ \ \ \ \ \ \ \ \ \ \ \ \ \ \ \ \ \ \ \ \ \ \ \  \=\mit
\footnotesize Takeo\ Ohsawa\\\mit
 \>\mit \footnotesize Mathematik  \>\mit \footnotesize Graduate School of Mathematics\\\mit
 \>\mit \footnotesize Universit\"at Wuppertal  \>\mit \footnotesize Nagoya University\\\mit
 \>\mit \footnotesize Gausstr. 20  \>\mit \footnotesize Chikusa-ku, Furocho\\\mit
 \>\mit \footnotesize D-42097 Wuppertal  \>\mit \footnotesize Nagoya 464-01\\\mit
 \>\mit \footnotesize Germany  \>\mit \footnotesize JAPAN\\\mit
  \end{tabbing}\vspace{-2\baselineskip}\mit
\end{document}